\documentclass[11pt]{article}
\usepackage{amsmath,amsthm,amssymb,epsf}
\usepackage[dvips]{graphicx}

\include{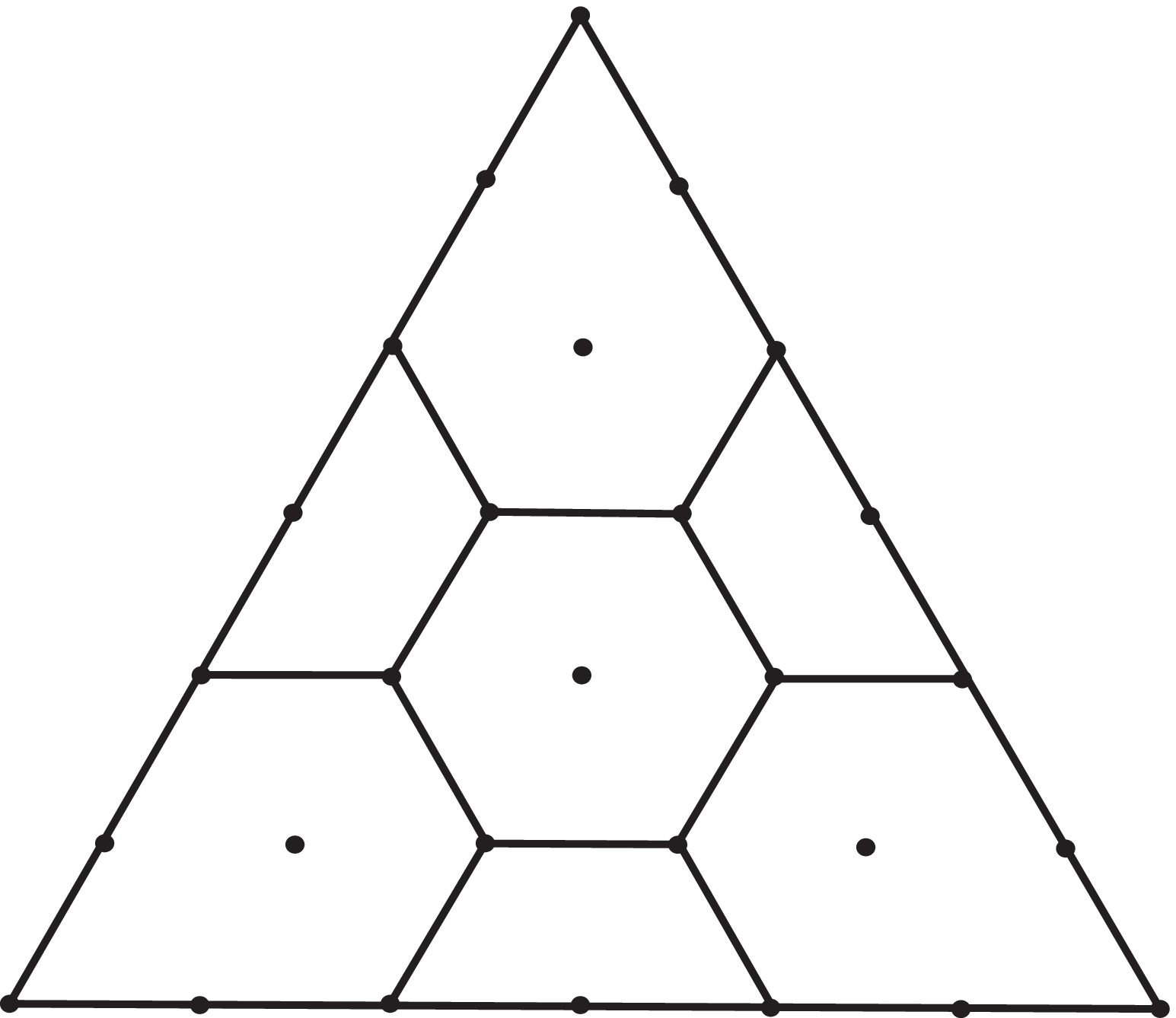}

\theoremstyle{plain}  

\newtheorem{thm}{Theorem}[section]
\newtheorem{prop}[thm]{Proposition}
\newtheorem{lem}[thm]{Lemma}   
\newtheorem{cor}[thm]{Corollary}

\theoremstyle{definition}

\newtheorem{defn}[thm]{Definition}

\theoremstyle{remark}

\newtheorem{rem}[thm]{Remark}

\DeclareMathOperator{\Spec}{Spec}
\DeclareMathOperator{\Proj}{Proj}
\DeclareMathOperator{\Def}{Def}

\DeclareMathOperator{\Hom}{Hom}

\DeclareMathOperator{\Ext}{Ext}

\DeclareMathOperator{\Pic}{Pic}

\DeclareMathOperator{\red}{red}

\DeclareMathOperator{\Hilb}{Hilb}

\DeclareMathOperator{\Lie}{Lie}

\DeclareMathOperator{\conv}{conv}
\DeclareMathOperator{\PGL}{PGL}

\DeclareMathOperator{\Cone}{Cone}

\DeclareMathOperator{\algebras}{algebras}
\DeclareMathOperator{\Sets}{Sets}
\DeclareMathOperator{\res}{res}
\DeclareMathOperator{\pr}{pr}

\newcommand{\QED}{\ifhmode\unskip\nobreak\fi\quad {\rm Q.E.D.}} 

\newcommand{\bA}{\mathbb A}
\newcommand{\bB}{\mathbb B}
\newcommand{\bC}{\mathbb C}

\newcommand{\bG}{\mathbb G}

\newcommand{\bN}{\mathbb N}
\newcommand{\bP}{\mathbb P}

\newcommand{\bR}{\mathbb R}
\newcommand{\bS}{\mathbb S}
\newcommand{\bT}{\mathbb T}
\newcommand{\bZ}{\mathbb Z}

\newcommand{\cA}{\mathcal A}
\newcommand{\cB}{\mathcal B}

\newcommand{\cI}{\mathcal I}
\newcommand{\cO}{\mathcal O}

\newcommand{\cN}{\mathcal N}

\newcommand{\cQ}{\mathcal Q}

\newcommand{\cS}{\mathcal S}
\newcommand{\cT}{\mathcal T}
\newcommand{\cU}{\mathcal U}
\newcommand{\cV}{\mathcal V}

\newcommand{\uP}{\underline{P}}

\newcommand{\cq}{/ \!  /}
\newcommand{\hq}{/ \! / \! /}
\newcommand{\omon}{M_{0,n}}
\newcommand{\mon}{\overline{M}_{0,n}}
\newcommand{\mgn}{\overline{M}_{g,n}}
\newcommand{\tbT}{\tilde{\mathbb T}}
\newcommand{\oo}{\overline{\Omega}}

\title{Compactification of the moduli space of hyperplane arrangements}
\author{Paul Hacking, Sean Keel and Jenia Tevelev}
\begin{document}
\maketitle

\section{Introduction}

Let $M^0$ denote the moduli space of arrangements of $n$ hyperplanes in $\bP^{r-1}_k$ in linear general position 
(i.e., ordered $n$-tuples of
hyperplanes in linear general position modulo the diagonal action of $\PGL(r)$). When $r =2$ the space, usually
denoted $\omon$, has a celebrated compactification due to Grothendieck and Knudsen, $\omon \subset \mon$, the moduli of
stable $n$-pointed rational curves.
The point of this note is to generalize the construction to higher dimensions. Of course $\mon$ is
the genus $0$ instance of $\mgn$, the moduli space of stable $n$-pointed curves of genus $g$. From
the point of view of Mori theory the correct 
generalisation of $\mgn$ is the moduli of semi log canonical pairs \cite{KSB88},\cite{Alexeev96a},\cite{Alexeev96b}, 
pairs $(S,B)$ of a variety with 
a boundary (a reduced Weil divisor) satisfying certain singularity assumptions generalizing toroidal
(we will not need the precise definition here). Such a space is expected to
exist in all dimensions, but known constructions depend on the minimal model program and so currently
apply only to varieties of dimension two or less. In this note we offer an alternative construction for
hyperplane arrangements (i.e., for generalizing $\mon$) which is quite elementary and which holds in
all dimensions. We will construct a projective scheme $M$, containing $M^0$ as open subset, 
and a flat projective family
$p:(\bS,\bB) \rightarrow M$ of (possibly reducible) $(r-1)$-dimensional 
varieties with boundary extending the universal family over $M^0$ (of ordered $n$-tuples of hyperplanes in 
$\bP^{r-1}$). The family has very nice properties:

\begin{thm} Let $(S,B = B_1 + \dots + B_n)$ be a fibre of $(\bS,\bB)$ over a closed point of $M$.
\begin{enumerate}
\item $(S,B)$ has at worst toroidal singularities (see Definition~\ref{toroidal}). The log canonical
sheaf $\omega_S(B)$ is a very ample line bundle, and the cohomology groups $H^i(S,\omega_S(B))$
vanish for $i > 0$. 
\item For each subset $I \subset \{1,2,3,\dots,n\}$ with $|I| = r-1$, the scheme-theoretic
intersection $\bB_I := \cap_{i \in I} \bB_I \subset \bS$ is a section of $p$, and the family
$(\bS,\bB)$ is semi-stable in a neighborhood of this section (i.e., near the corresponding point of
the fibre, $S$ and the $B_i$ are smooth and $B_1 + \dots + B_n$ has normal crossings).
\item The map given by taking residues along the sections $\bB_I$ 
$$
\res: p_*\omega_p(\bB) \rightarrow \oplus_{I} p_*\cO_{\bB_I} = \wedge^{r-1} k^n \otimes \cO_M
$$
is an isomorphism onto $\wedge^{r-1}h^* \otimes \cO_M \subset \wedge^{r-1} k^n \otimes \cO_M$, 
where $h=k^n/ (k \cdot (1,\dots,1))$. In particular $p_*\omega_p(\bB)$ is locally free of rank
${n-1}\choose {r-1}$. Its formation commutes with all base-extensions. In particular the above
residue map determines a basis of $H^0(S,\omega_S(B))$ canonically associated
to the pair $(S,B)$. 
\item The global sections given by $\res$ induce a 
canonical embedding 
$$
\bS \subset M \times G(r,n) \subset M \times \bP(\wedge^{r} k^n)
$$
where $G(r,n) \subset \bP(\wedge^{r} k^n)$ is the Pl\"ucker embedding of the Grassmannian
of $r$-planes in affine $n$-space. The induced map
$M \rightarrow \Hilb(G(r,n))$ is a closed embedding. The closure of $M^0 \subset M$ is 
Kapranov's Chow quotient $G(r,n) \cq H$ of the Grassmannian by its maximal torus, see \cite{Kapranov93}.
\item $p$ is a flat family of log canonically polarised semi log canonical pairs and so
defines a map from $M$ to the moduli stack of
semi log canonical pairs. This is a closed immersion.
\end{enumerate}
\end{thm}
Furthermore, the family $(\bS,\bB)$ is universal and identifies $M$ as a natural moduli 
space of pairs satisfying properties as in the theorem --- what we call `very stable pairs'.
See Section~\ref{verystable} for the precise statement.

Unfortunately our $M$ will not in general be irreducible, see Section~\ref{excpt}, and thus is not precisely a 
compactification of $M^0$. We do not know a functorial characterisation of the closure of $M^0$ (i.e., of the Chow
quotient $G(r,n) \cq H$). 

\subsection{General Philosophy}

Before turning to the technical details let us outline the general idea, which is adapted from
ideas of \cite{Kapranov93} and \cite{Lafforgue03}. Begin first with a pair
$(S,B= B_1 + \dots B_n)$ of $\bP^{r-1}$ together with $n$ hyperplanes in linear general position. 
The main observation is that moduli of
such pairs can be identified with moduli of equivariant embeddings of a fixed toric variety ---
the normal projective toric variety associated to the so called hypersimplex $\Delta(r,n)$ --- 
in the Grassmannian, $G(r,n)$.

By the Gel'fand-MacPherson transform $M^0$ is identified with the set 
of orbits $G^0(r,n)/H$, where $G^0 \subset G(r,n)$ is the open subset where all Pl\"ucker
coordinates are non-zero and $H=\bG_m^n/\bG_m \subset \PGL(n)$ is the standard maximal torus. 
In \cite{Kapranov93} this correspondence is formulated elegantly as follows: 
A choice of
linear equations for the hyperplanes yields an embedding 
$S \subset \bP^{n-1}$ so that the configuration $B$ is the restriction of the coordinate hyperplanes.
$H$ acts freely on the orbit of $[S] \in G(r,n)$ so we have an isomorphism
$$
m: H \rightarrow H \cdot [S], \, h \mapsto h^{-1} [S].
$$
Observe $S \setminus B \subset H$ is identified with
$$
\{P \in H \cdot [S] \ | \ e \in P\} = H \cdot [S] \cap G(r-1,n-1)_e
$$
where $e = (1,\dots,1) \in \bP^{n-1}$, and $G(r-1,n-1)_e \subset G(r,n)$ is
the sub Grassmannian of $r$-planes that contain the fixed vector $e$. This identification
is easily seen to extend to the closure (and indeed to degenerations), see Section~\ref{univfamily}, 
$S = \overline{H \cdot [S]} \cap G(r-1,n-1)_e$. This realizes $S$ as 
a complete intersection inside the orbit closure $\overline{H \cdot [S]}$, the normal
projective toric variety corresponding to the polytope $\Delta(r,n)$. 
Kapranov calls the orbit closure a \emph{Lie complex} and $S \subset \overline{H \cdot [S]}$ 
its \emph{visible contour}. This realizes $M^0$ as a locus in $\Hilb(G(r,n))$ of
generic orbit closures. The closure of this locus is Kapranov's Chow quotient compactification 
$M^0 \subset G(r,n) \cq H$. By definition $G(r,n) \cq H$ 
carries a flat family, with generic fibre these orbit closures. The advantage of the
approach is that the degenerate fibres are quite easy to understand --- the generic fibres are
closures of generic $H$-orbits and 
are embeddings of the normal projective toric variety associated to $\Delta(r,n)$, 
special fibres are reduced unions of (top dimensional) orbit closures, which
are normal projective toric
varieties associated to cells in certain tilings (called matroid decompositions) of $\Delta(r,n)$,
see Corollary~\ref{toricfamily}. In particular we have a flat family of pairs $(\bT,\bB_T)$ of \emph{broken} toric varieties
and their toric boundaries. 
A simple but clever observation of Lafforgue shows that the visible contour construction 
extends to all of $(\bT,\bB_T)$ --- and yields exactly as above a flat family 
$(\bS,\bB) \subset (\bT,\bB_T)$ of complete intersections, transverse to the 
toric boundary, and in particular $(\bS,\bB)$ a flat family of pairs 
with toroidal singularities, compactifying the universal family
of hyperplane arrangements over $M^0$. See Section~\ref{Laftrans}. We observe that for each fibre
$(S,B)$ of  $(\bS,\bB)$ the Pl\"ucker embedding
$$
S \subset G(r-1,n-1)_e \subset G(r,n) \subset \bP(\wedge^r k^n)
$$
(and so the Hilbert point $[S] \in \Hilb(G(r,n))$) is given by a canonical
basis of global log canonical forms, and in particular is canonically determined by
the isomorphism class of the pair $(S,B)$, see Theorem~\ref{residues}. In this way
$(\bS,\bB) \rightarrow G(r,n) \cq H$ induces a closed immersion of 
$G(r,n) \cq H$ into the moduli stack of semi log canonical pairs; thus $G(r,n) \cq H$ 
is a sub moduli space of pairs. Unfortunately we cannot identify
the image --- we do not know precisely which semi log canonical pairs are limits of 
generic hyperplane arrangements. 
Here we use an alternative construction: 
Instead of $G(r,n)\cq H$ we make use of $M \subset \Hilb(G(r,n))$, a closed 
subscheme of the so called toric Hilbert scheme, see \cite{HS04}. $M$ parameterises $\bG_m^n$-equivariant
closed subschemes  of $\tilde{G}(r,n)$ (the cone over the Grassmannian in its Pl\"ucker embedding) 
with a prescribed multigraded Hilbert function, see Section~2. 
$M^0$ immerses in $M$ as an
open subset, with closure $G(r,n) \cq H$, and, because the toric Hilbert scheme represents
a natural functor, $M$ admits a functorial description as a moduli space of pairs with toroidal
singularities (satisfying various other properties), which we call very stable pairs. 
See Section~\ref{verystable} for the precise statement. 

\section{The log canonical model of the complement of a hyperplane arrangement} \label{lcmodel}

This short section is not logically required for the proof of the main theorem --- everything we do here
we'll redo in later sections in greater generality. As we think the construction is of independent interest,
we have written the section so that it can be read on its own, at the cost of some subsequent repetition. 

We describe an explicit compactification $(S,B)$ of the complement $U$ of a 
hyperplane arrangement, following \cite{Kapranov93}. We show that $(S,B)$ is the log canonical model of 
$U$, i.e., the canonical compactification of the algebraic variety $U$ obtained via the minimal model 
program. These compactifications occur as the components of the fibres of the universal family 
$(\bS,\bB)/M$. 

Let $\cA=(H_1,\ldots,H_n)$ be an (ordered) arrangement of hyperplanes in 
$\bP^{r-1}$. Let $U=\bP^{r-1} \setminus \cup \cA$, the complement.
Assume that the stabiliser of $\cA$ in $\PGL(r)$ is finite.
Equivalently, the matroid of $\cA$ is \emph{connected} \cite{GS87}, i.e., 
there does not exist a decomposition $k^r=V_1 \oplus V_2$ such that 
for each $i$ either $\bP(V_1) \subset H_i$ or $\bP(V_2) \subset H_i$.

Choose homogeneous equations $F_i$ for the $H_i$, and consider the linear 
embedding
$$F=(F_1:\ldots:F_n) : \bP^{r-1} \subset \bP^{n-1}.$$
Let $H=\bG_m^n/\bG_m \subset \bP^{n-1}$ be the usual torus embedding.
Observe that the embedding $F$ is determined up to translation by an element 
of $H$, and restricts to a (closed) embedding $U \subset H$.

Let $G(r,n)$ denote the Grassmannian of $r$-planes in $k^n$.
Let $V$ denote the $H$-orbit in $G(r,n)$ determined by $F$.
The weight polytope of $V$ is the matroid polytope of $\cA$.
It has full dimension $n-1$ since by assumption the matroid of $\cA$ is 
connected 
(see \cite{GS87}), and its vertices affinely generate the lattice 
(see, e.g., \cite{Kapranov93}, p.~47, Proof of Prop.~1.2.15).
Hence $H$ acts freely on $V$. The embedding $U \subset V$ given by 
$$u \mapsto F(u)^{-1}[F(\bP^{r-1})]$$
is canonical (it does not depend on the choice of $F$). 

Let $G(r-1,n-1)_e \subset G(r,n)$ denote the locus of subspaces containing the 
vector $e=(1,\ldots,1) \in k^n$. Note that $G(r-1,n-1)_e$ is identified with 
the Grassmannian of $(r-1)$-planes in $h = k^n / k \cdot e$, the Lie algebra 
of $H$.
Observe that the locus $U \subset V$ in $G(r,n)$ equals $V \cap G(r-1,n-1)_e$.

Let $S$ and $T$ denote the closures of $U$ and $V$ in $G(r,n)$, respectively.
The variety $T$ is isomorphic to the normal toric variety associated to the 
matroid polytope of $\cA$. Write $B=S \setminus U$ and $B_T=T \setminus V$, the toric 
boundary of $T$.

\begin{lem}[Lafforgue, cf. Thm.~\ref{transthm}] $S$ is equal to the scheme-theoretic intersection $T \cap G(r-1,n-1)_e$.
The multiplication map $H \times S \rightarrow T$ is smooth.
\end{lem}
\begin{proof} Let $S'=T \cap G(r-1,n-1)_e$, then $S$ is clearly an irreducible component of $S'$. 
We show that $S'$ is reduced and irreducible, so $S=S'$.
Let $\bP(\cU) \subset  G(r,n) \times \bP^{n-1}$ denote the projectivised universal bundle over $G(r,n)$.
The multiplication map $H \times S' \rightarrow T$ is identified with the projection 
$\bP(\cU) \cap T \times H \rightarrow T$ 
via $(h,s) \mapsto (hs,h).$
In particular $H \times S' \rightarrow T$ is smooth, and $S'$ is reduced and irreducible.
\end{proof}

\begin{thm} \label{lcmodelthm}
$(S,B)$ is the log canonical model of $U$. Moreover,
\begin{enumerate}
\item $(S,B)$ has toric singularities.
\item $K_S+B$ is very ample. 
\item The embedding $S \subset G(r-1,n-1)_e$ is given by the locally free sheaf
$\Omega_S(\log B)$ and the map
$$h^* \rightarrow H^0(\Omega_S(\log B)), \, 
(\lambda_1,\ldots,\lambda_n) \mapsto \sum \lambda_i \frac{dF_i}{F_i}.$$
\end{enumerate}
\end{thm}
\begin{proof}
$(S,B)$ has toric singularities by the Lemma.
Assuming (3), $\Omega_S(\log B)$ is identified with the restriction of the 
dual of the universal sub-bundle $\cU_e \subset \cO_{G_e} \otimes h$ on 
$G(r-1,n-1)_e$. So $\omega_S(B)=\wedge^{r-1}\Omega_S(\log B)$ is identified 
with the restriction of the Pl\"{u}cker line bundle on $G_e$. 
Hence $K_S+B$ is very ample.

For $P \in H$, let $\mu_P : H \rightarrow H$ be the map given by 
multiplication by $P$.
The embedding $U \subset G(r-1,h)$
is the Gauss map associated to the embedding $U \subset H$, i.e., the map
$$g : U \rightarrow G(r-1,h), P \mapsto [d(\mu_P^{-1})T_P U].$$
Indeed, since $U \subset H$ is the restriction of the \emph{linear} embedding
$\bP^{r-1} \subset \bP^{n-1}$, all the tangent spaces $T_P U$ are equal to 
$\bP^{r-1} \subset \bP^{n-1}$ (when regarded as subspaces of $\bP^{n-1}$).
An explicit computation shows that the embedding $U \subset G(r-1,h)$ is given
by the surjection
$$
h^* \otimes \cO_U \rightarrow \Omega_U, \,
(\lambda_1,\ldots,\lambda_n) \mapsto \sum \lambda_i \frac{dF_i}{F_i}.
$$
This map extends to the surjection
$$h^* \otimes \cO_S = \Omega_T(\log B_T)|_S \rightarrow \Omega_S(\log B).$$
given by the embedding $S \subset T$. Statement (3) follows.
\end{proof}

If $k=\bC$,  part (3) may be explained conceptually as follows.
The exponential map 
$$\exp: h \rightarrow H, \, (\lambda_1,\ldots,\lambda_n) \mapsto (\exp(\lambda_1), \ldots \exp(\lambda_n))$$
identifies $H$ with the quotient $h/(2 \pi i)N$, where $N =\bZ^n/\bZ e \subset h =\bC^n/\bC e$, 
the cocharacters of $H$.
Assume for simplicity that the hyperplanes $H_1,\ldots,H_n$ are distinct, then the map 
$h^* \rightarrow H^0(\Omega_S(\log B))$ is an isomorphism.
The embedding $U \subset H$ is identified with the (generalised) Albanese map
$$U \rightarrow H^0(\Omega_S(\log B))^*/H_1(U,\bZ), \, 
P \mapsto \left( \omega \mapsto \int^P_{P_0} \omega \right),$$
where $P_0 \in U$ is a fixed basepoint. 
Recall that $g: U \rightarrow G(r-1,h)$ is the Gauss map for $U \subset H$.
Using the integral formula for the embedding $U \subset H$ and the fundamental
theorem of calculus, we deduce that $U \subset G(r-1,h)$ is given by the locally free sheaf 
$\Omega_U$ and the surjection
$$h^* \otimes \cO_U = H^0(\Omega_S(\log B)) \otimes \cO_U \rightarrow \Omega_U.$$
The result follows as above.

\section{Construction of the moduli space of pairs}

\subsection{Multigraded Hilbert schemes}
$M$ is a multigraded Hilbert scheme as defined in \cite{HS04}. 
We briefly review the definition and basic properties.

Let $T=\oplus_{a \in A} T_a$ be a $k$-algebra graded by an Abelian group $A$.
Fix a function $h : A \rightarrow \bN$. For $R$ a $k$-algebra, let $H^h_T(R)$
be the set of $A$-homogeneous ideals $I \subset T \otimes R$ such that, for each $a \in A$,
$T_a \otimes R/ I_a$ is a locally free $R$-module of rank $h(a)$. This defines a functor 
$H^h_T: (k-\algebras) \rightarrow (\Sets)$. It is represented by a quasiprojective scheme over $k$, the 
\emph{multigraded Hilbert scheme} $H^h_T$.
If $T$ is a polynomial ring and the multigrading is positive (i.e., $T_0=k$), then $H^h_T$ is projective.

Let $S=k[x_1,..,x_N]$ and $\bA=\Spec S$. Fix a map 
$$\phi : \bZ^N \rightarrow \bZ^n, \ e_i \mapsto a_i$$
corresponding to a homomorphism of tori $\bG^n_m \rightarrow \bG^N_m$, where $\bG^N_m$ is the big torus acting on $\bA$.
Let $\cA=\{a_1,\ldots,a_N\}$, the set of weights for the torus action $\bG_m^n \curvearrowright \bA$, 
and $A=\bZ\cA \subset \bZ^n$ the lattice generated by $\cA$.
The map $\phi$ defines an $A$-grading of $S$ such that the $A$-homogeneous ideals $I \subset S$ are the ideals defining 
$\bG^n_m$-invariant closed subschemes in $\bA$.

Let $\bN \cA \subset A$ be the semigroup generated by $\cA$.
Define $h : A \rightarrow \bN$ by $h(a) = 1$ if $a \in \bN\cA$  and $h(a)=0$ otherwise.
The multigraded Hilbert scheme $H^h_S$ is the {\it toric Hilbert scheme} for the torus action 
$\bG_m^n \curvearrowright \bA$ \cite[Sec.~5]{HS04}. Roughly speaking, $H^h_S$ parameterises generic $\bG^n_m$-orbit closures in $\bA$ 
and their toric degenerations. More precisely, let $X_{\cA}$ denote the orbit closure 
$\overline{\bG^n_m  \cdot e} \subset \bA$, where $e=(1,\ldots,1) \in \bA$. 
Then $X_{\cA}$ defines a distinguished point $[X_{\cA}] \in H^h_S$, and the orbit closure 
$\overline{\bG^N_m \cdot [X_{\cA}]} \subset H^h_S$ is an irreducible component of $H^h_S$.

If $X=\Spec T \subset \bA$ is a $\bG^n_m$-invariant closed 
subscheme, then $T$ is $A$-graded and $H^h_T$ is the closed subscheme of $H^h_S$ parameterising subschemes of $X$.

\subsection{Stable toric varieties}

A subscheme $Z \subset \bA$ defining a point of the toric Hilbert scheme $H^h_S$ is an affine stable toric variety 
as defined in \cite{Alexeev02} (assuming $Z$ is seminormal and reduced and the multigrading is positive). 
We review the construction of stable toric varieties.

Let $A$ be a lattice and $\Omega$ a subdivision of a rational polyhedral cone $\omega$ in $A_{\bR}$.
For $\sigma \in \Omega$ let $R_{\sigma}$ denote the semigroup algebra $k[\sigma \cap A]$ and 
$T_{\sigma} \subset X_{\sigma}$ the associated torus embedding. 
Fix glueing data $t_{\sigma\tau} \in T_{\tau}$ for each $\tau \subset \sigma$ 
satisfying the compatibility
condition $t_{\tau\upsilon} \cdot t_{\sigma\tau} = t_{\sigma\upsilon}$ in $T_{\upsilon}$
for each triple $\upsilon \subset \tau \subset \sigma$.
Define $p_{\sigma \tau}=t_{\sigma \tau}\circ\pr_{\sigma \tau}$ for $\tau \subset \sigma$,
where $\pr_{\sigma \tau}$ is the canonical surjection $R_{\sigma} \rightarrow R_{\tau}$.
Finally, let $R[\Omega,t]$ be the inverse limit of the system $(R_{\sigma}, p_{\sigma \tau})$.

\begin{rem}
Equivalently, $R[\Omega,t]$ is the equaliser of the maps $\oplus R_{\sigma} \rightrightarrows \oplus R_{\tau}$,
where the direct sums are over maximal cones $\sigma \in \Omega$ and codimension $1$ interior cones $\tau \in \Omega$, 
respectively. That is, $R[\Omega,t]$ is the subalgebra of $\oplus R_{\sigma}$ consisting of elements $f=(f_{\sigma})$ 
such that $p_{\sigma_1\tau}(f_{\sigma_1})=p_{\sigma_2\tau}(f_{\sigma_2})$ for each pair $\sigma_1,\sigma_2$ of maximal cones 
meeting in a common facet $\tau$.
\end{rem}

The variety $X=X(\Omega,t):=\Spec R[\Omega,t]$ has irreducible components $X_{\sigma}=\Spec R_{\sigma}$ for 
$\sigma \in \Omega$
a maximal cone. Combinatorially, the $X_{\sigma}$ are glued to form $X$ in the same way that the cones $\sigma$ are 
glued to form $\omega$. 
That is, for each maximal cone $\sigma$, the facets of the cone $\sigma$ correspond to the irreducible components
of the toric boundary $X_{\sigma} \backslash T_{\sigma}$ of $X_{\sigma}$, and 
if $\sigma_1$ and $\sigma_2$ meet in a common facet then $X_{\sigma_1}$ and $X_{\sigma_2}$ are glued along the corresponding 
divisor. Note that there are also continuous glueing parameters determined by $t$.
There is an action of the torus $T=\Hom(A,\bG_m)$ on $X$ extending the action on each component.
The algebra $R[\Omega,t]$ with its corresponding $A$-grading has Hilbert function $h(a)=1$ for 
$a \in \omega \cap A$ and $h(a)=0$ otherwise.

\begin{defn} An affine stable toric variety is a variety with torus action 
of the form $T \curvearrowright X(\Omega,t)$ for some $\Omega,t$.
\end{defn}

\begin{rem} If $t_{\sigma \tau}=1$ for each $\tau \subset \sigma$, then $R[\Omega,t]$ can be 
alternatively described as follows, cf. \cite{Stanley87}. 
As a $k$-vector space, $R=\oplus \ k \cdot \chi^a$ where the sum is over the semigroup 
$\omega \cap A$. The ring structure on $R$ is defined by 
$\chi^a \cdot \chi^b = \chi^{a+b}$ if $a$ and $b$ are contained in some cone $\sigma \in \Omega$, and 
$\chi^a \cdot \chi^b=0$ otherwise. 
\end{rem}

Let $M$ be a lattice, $P \subset M_{\bR}$ a polytope with integral vertices, and $\uP$ a subdivision of $P$.
Let $A=M \oplus \bZ$, and embed $P$ in the affine hyperplane $M_{\bR} \oplus 1 \subset A_{\bR}$. 
Let $\Omega$ be the fan of cones over faces of $\uP$.
Fix glueing data $t$ as above and define $Y=Y(\uP,t):=\Proj R[\Omega,t]$.
The irreducible components of $Y$ are the polarised projective toric varieties 
$Y_{P'}=\Proj R_{\Cone(P')}$ associated to the maximal polytopes $P' \in \uP$.
The combinatorics of the glueing of the $Y_{P'}$ is encoded by $\uP$.
There is an action of the torus $H=\Hom(M,\bG_m)$ on $Y$, and the polarisation $\cO(1)$ on $Y$ has a natural 
$H$-linearisation.

\begin{defn} A polarised stable toric variety is a projective variety with a torus action together
with a linearised ample sheaf of the form $H \curvearrowright (Y(\uP,t), \cO(1))$
\end{defn}

\begin{rem}
In \cite{Alexeev02} the definition of stable toric varieties is more general, 
and the special case above is referred to as the ``convex 1-sheeted case''.
\end{rem}

\subsection{The construction}

Let $G(r,n) \subset \bP=\bP(\wedge^r k^n)$ be the Pl\"ucker embedding of the Grassmannian of $r$-planes in $k^n$.
Let $\tilde{G}(r,n) \subset \bA$ be the cone over the Pl\"ucker embedding, and
$S$ and $T$ the coordinate rings of $\bA$ and $\tilde{G}(r,n)$ respectively. 
Let $\bG_m^n \curvearrowright \bA$ be the standard $\bG^n_m$-action and $H^h_S$ the associated 
toric Hilbert scheme.
                                 
\begin{defn} Let $M=H^h_T$, the closed subscheme of 
the toric Hilbert scheme $H^h_S$ parametrising subschemes of $\tilde{G}(r,n)$.
\end{defn}
Note immediately that we have an open immersion $M^0 \subset M$ given by the Gel'fand--MacPherson correspondence 
$M^0 =  G^0(r,n)/H$.

The set of weights of $\bG^n_m \curvearrowright \bA$ is  
$$\cA=\left\{ e_{i_1}+\cdots+e_{i_r} \ \big| \ i_1 < \cdots < i_r \right\} \subset \bZ^n $$
where $e_1,\ldots,e_n$ is the standard basis of $\bZ^n$. 
The set $\cA$ is the set of vertices of the \emph{hypersimplex}
$$
\Delta(r,n) :=\left\{(x_1,\dots,x_n) \in \bR^n\ \big|\ \sum x_i = r,\ 0 \leq x_i \leq 1\right\}. 
$$
The polytope $\Delta(r,n)$ has $2n$ facets $(x_i=0)$ and $(x_i=1)$, $i=1,\ldots,n$. 
Write $P=\Delta(r,n)$.

We consider polytopes $P' \subset P$ which are the convex hull of a subset of the vertices $\cA$ of $P$.
We regard the coordinates of $\bA$ as labelled by $\cA$.
For $P' \subset P$, let $x_{P'} \in \bA$ be the point with coordinates $1$ for $a \in P' \cap \cA$ and 
$0$ otherwise, 
and $X_{P'}$ the orbit closure $\overline{\bG^n_m \cdot x_{P'}}$. $X_{P'}$ is the affine toric variety 
(possibly non-normal) associated to the semigroup $\bN(P' \cap \cA) \subset A$ generated by $P' \cap \cA$. 

Let $\tbT \subset \tilde{G}(r,n) \times M$ denote the universal family over $M$.

\begin{thm} \label{stvfamily}
Each fibre of $\tbT/M$ is a reduced affine stable toric variety associated to a subdivision of $\Cone(P)$ 
induced by a subdivision of $P$ into matroid polytopes.
\end{thm}

\begin{proof} 
Let $Z$ be a fibre of $\tbT/M$.
By \cite[10.10]{Sturmfels95} there is a polyhedral subdivision $\uP$ of $P$ such that
$\red{Z}= \bigcup_{P'} Z_{P'}$ where the union is over maximal polytopes $P' \in \uP$, and $Z_{P'}$ is a translate of 
$X_{P'}$ by the big torus acting on $\bA$.

Each $P'$ is a matroid polytope since $Z \subset \tilde{G}(r,n)$.
Hence the set $P' \cap \cA$ generates the saturated semigroup $\Cone(P') \cap A$ by \cite{White77}, 
so $Z_{P'}$ is normal. 
It also follows that $Z$ is reduced. For, we have the surjections of coordinate rings
$$k[Z] \rightarrow k[\red Z] \rightarrow k[Z_{P'}]$$ 
and $\dim_k k[Z]_a = \dim_k k[Z_{P'}]_a=1$ for $a \in \Cone(P') \cap A$. 
Thus $k[Z]_a=k[\red Z]_a$ for each $a \in A$ and $Z=\red Z$ as claimed.

If $P'_1$ and $P'_2$ intersect in a common facet, the corresponding boundary divisors of 
$Z_{P'_1}$ and $Z_{P'_2}$ coincide with the scheme-theoretic intersection $Z_{P'_1} \cap Z_{P'_2}$.
Indeed, the ideal of $Z_{P'} \subset Z$ is the direct sum of the graded pieces $k[Z]_a$ of $k[Z]$ 
for $a \notin \Cone{P'}$. We deduce that $k[Z]$ is the equaliser of the maps
$$\oplus k[Z_{P'}] \rightrightarrows \oplus k[Z_{P''}],$$
where the $Z_{P''}$ are the strata of $Z$ corresponding to interior codimension 1 faces $P'' \in \uP$.
Hence $Z$ is an affine stable toric variety. 
\end{proof}

\begin{cor} The natural map 
$M \rightarrow \Hilb(G(r,n))$ obtained by projectivising $\tilde{G}(r,n) \subset \bA$ is a closed embedding.
\end{cor}
\begin{proof}
Let $Z \subset \tilde{G}(r,n) \times \Spec R$ be an $R$-valued point of $H^h_T$. The family $Z/R$ is flat and has reduced 
fibres by Theorem~\ref{stvfamily}.
It follows by \cite[2.32]{Matsumura89} that the ideal $I \subset S \otimes R$ of $Z \subset \bA \times \Spec R$ is 
saturated. Hence the map $H^h_T \rightarrow \Hilb(G(r,n))$ is an injection on $R$-points for each $R$.
\end{proof}

\begin{cor} The closure of $M^0 \subset M$ is the Chow quotient $G(r,n) \cq H$. 
\end{cor}
\begin{proof} 
By definition $G(r,n) \cq H$ is the closure of $M^0$ in $\Hilb(G(r,n))$.
\end{proof}

Let $\bT \subset G(r,n) \times M$ denote the family obtained by projectivising $\tbT \subset \tilde{G}(r,n) \times M$.

\begin{cor} \label{toricfamily}
Each fibre of $\bT/M$ is a reduced projective stable toric variety associated to a subdivision of $P$ into 
matroid polytopes.
\end{cor}

\subsection{Relation to Lafforgue's space} \label{Lafspace}

Lafforgue defines a projective scheme $\oo=\oo^{\Delta(r,n)}$ with an open immersion $M^0 \subset \oo$. 
It may be constructed as follows (see \cite[2.9]{KT04}). Let $\bP \hq_n H \rightarrow \bP \hq H \subset \Hilb(\bP)$ 
be the normalisation of the \emph{Hilbert quotient} of $\bP(\wedge^r k^n)$, i.e., the closure in $\Hilb(\bP)$ of the 
locus of generic $H$-orbit closures. The space $\oo$ is the inverse image in $\bP \hq_n H$ of 
$\bP \hq H \cap \Hilb(G(r,n))$. This construction induces a finite map $\oo \rightarrow M$ such that the family over 
$\oo$ (coming from $\Hilb(\bP)$) is the pullback of $\bT$. It is an isomorphism over $M^0 \subset M$.

Roughly speaking, the space $\oo$ is a moduli space of varieties with log structures --- see \cite[Ch.~5]{Lafforgue03}
for the precise statement. 
Our space $M$ is a sub moduli space of stable pairs, see Section~\ref{verystable}.
Given a $k$-point $[(S,B)] \in M$, a point of $\oo$ over $[(S,B)]$ is given by a smooth log structure on $S/k$ 
which is nontrivial over the divisors $B_i \subset S$ and the singular locus. 
Such log structures do not always exist, see Section~\ref{excpt}.
Moreover, we expect that the log structure is not unique in general, i.e., the map $\oo \rightarrow M$ is not injective on $k$-points.

The space $\oo$ is in general reducible by \cite[3.13]{KT04}, which of course implies the same for $M$. 
$M$ also has components outside (the image) of $\oo$, see Section~\ref{excpt}. Ideally, we would like $M$ to be a 
connected component of the moduli space of stable pairs, but we do not know if this is the case.

\section{Construction of universal family of pairs} \label{univfamily}

\subsection{Lafforgue Transversality}\label{Laftrans}

Section~4.1 and Theorem~\ref{transthm} are based on \cite[5.1]{Lafforgue03}. 
Let $G$ be a scheme on which an algebraic group $W$ acts.
Let $\cV \subset G \times W$ be a closed $W$-equivariant 
subscheme. Define $\cV_e := \cV \cap (G \times \{e\})$, where $e \in W$ is the identity element. 
Note the first projection $\cV_e \rightarrow G$ is a closed embedding. 
Let $G_{e,\cV} \subset G$ be the image. 
\begin{lem} \label{Laftranslem}
The multiplication map
$\cV_e \times W \rightarrow \cV$ is an isomorphism, and identifies 
the multiplication map $G_{e,\cV} \times W \rightarrow G$ with the
first projection $\cV \rightarrow G$. 

Let $G' \rightarrow G$ be an $W$-equivariant map, and let
$\cV' \subset G' \times W$ be the pullback. Then 
$G'_{e,\cV'} \subset G'$ in the pullback of $G_{e,\cV} \subset G$.
\end{lem}
\begin{proof} The map $\cV \rightarrow \cV_e \times W$ given
by $(g,w) \rightarrow ((gw^{-1},e),w)$ is easily seen to be
inverse to right multiplication. The rest is easy to check.
\end{proof}

\begin{rem} Of course if $\cV \rightarrow G$ is smooth, then by the Lemma so
is the map $G_{e,\cV} \times W \rightarrow G$. 
\end{rem}

\subsection{Visible contours}

Now let $G_e=G(r-1,n-1)_e \subset G=G(r,n)$ be the locus of subspaces containing $e=(1,\ldots,1)$.
Let $H= \bG_m^n / \bG_m$ be the standard maximal torus in $\PGL(n)$ and $h=k^n/ k \cdot e$ the Lie 
algebra of $H$.
Note that $G_e$ is identified with $G(r-1,h)$.

\begin{defn} 
Following \cite{Kapranov93}, we define the \emph{family of visible contours} $p : (\bS,\bB) \rightarrow M$ as follows.
Let $\bS$ denote the scheme-theoretic intersection $\bT \cap (G_e \times M)$.
Let $\bB_T$ denote the relative toric boundary of $\bT/M$ and $\bB$ its restriction
to $\bS$. 
\end{defn}

There is a decomposition $\bB_T= \sum_{i=1}^n \bB^+_{i,T} + \sum_{i=1}^n \bB^-_{i,T}$,
where $\bB^+_{i,T}$ and $\bB^-_{i,T}$ are the components of the $\bB_T$ corresponding to the facets $(x_i=1)$ and 
$(x_i=0)$ of $\Delta(r,n)$ respectively. The components $\bB^-_{i,T}$ are disjoint from $\bS$, so $\bB = \sum_{i=1}^n \bB_i$ where 
$\bB_i:=\bB^+_{i,T}|_{\bS}$. 

The family $(\bS,\bB_1+\cdots+\bB_n)/M$ extends the universal family of hyperplane arrangements over $M^0$ by \cite[3.2.3]{Kapranov93} 
or Section~\ref{lcmodel}.

\begin{defn} \label{toroidal} 
A \emph{toroidal pair} $(S,B)$ is a (possibly reducible) variety $S$ together with a reduced divisor $B \subset S$ which
is \'{e}tale locally isomorphic to a stable toric variety with its toric boundary.
\end{defn}

\begin{thm}\label{transthm} The multiplication map
$\bS \times H \rightarrow \bT$ is smooth with image $\bT \setminus \cup \bB_{i,T}^-$. 
The family $\bS$ and the $\bB_i$ are flat over $M$. 
The embedding $\bS \subset \bT$ is the pullback of $G(r-1,n-1)_e \subset G(r,n)$. It is a regular embedding
with normal bundle the restriction of the universal rank $n-r$ quotient bundle on $G(r-1,n-1)_e$.
\end{thm}
\begin{proof} Let $\cU \subset G(r,n) \times k^n$ be the universal rank $r$ bundle,
and $\cV \subset \cU$ the intersection of $\cU$ with the diagonal torus $W \subset k^n$ (the locus
where all coordinates are non-zero). Then by definition
$\bS \subset \bT$ is the pullback of $G(r-1,n-1)_e \subset G(r,n)$ and, following the
notation of Section~\ref{Laftrans}, $G_e = G_{e,\cV}$. Now it follows from Lemma~\ref{Laftranslem} that 
$\bS \times W \rightarrow \bT$ is identified with the pullback
of $\cV \rightarrow G(r,n)$, and in particular is smooth. Since the
scalar matrices act trivially on $G(r,n)$, and thus on $\bS$, it
follows that $\bS \times H \rightarrow \bT$ is smooth as well. A particular closed stratum of a fibre
$T$ of $\bT$ is disjoint from the visible contour (or equivalently the image of 
$\bS \times H$) iff the affine $r$-plane corresponding to its generic point lies in a coordinate hyperplane
of $k^n$, which holds iff the corresponding matroid polytope lies in a face $(x_i =0)$ of $\Delta(r,n)$
and thus
iff the stratum lies in the boundary divisor $\bB_{i,T}^{-}$ (see, e.g., \cite[Sec.~1]{Kapranov93}). 
Since $H$ acts trivially
on $M$, flatness of $\bS$ and the $\bB_i$ (over $M$) now follow from flatness of
$\bT$ and the components of $\bB_T$. 

Finally, the closed subscheme $G_e \subset G$ is the zero locus of the section $\bar{e}$ of the quotient bundle $\cQ$ given by $e \in k^n$,
thus $G_e \subset G$ is a local complete intersection with normal bundle $\cN_{G_e/G}=\cQ|_{G_e}=\cQ_e$.
Now by the previous results $\bS \subset \bT$ 
is also a local complete intersection with normal bundle
$\cN_{\bS/\bT}=\cQ_e|_{\bS}$.
\end{proof}

\begin{cor}
Let $(T,B_T)$ be a fibre of $(\bT,\bB_T)/M$ and $(S,B)$ its visible contour.
\begin{enumerate}
\item $(S,B)$ has toroidal singularities.
\item Consider the stratification of $S$ induced by the stratification of $T$ by orbit closures.
A stratum $S'=S \cap T'$ is non-empty if and only if $T' \not\subset \bigcup B_{i,T}^-$.
In this case, $S'$ is irreducible and normal of the expected dimension $\dim T' - (n-r)$.
\end{enumerate}
\end{cor}

\begin{rem}
The stratification of $S$ coincides with that defined by arbitrary intersections of components of $S$ and $B$. 
In particular, it is obviously intrinsic. 
Let $\underline{P}$ be the polyhedral subdivision of $P=\Delta(r,n)$ associated to the stable toric variety $T$.
The poset of orbit closures in $T$ is identified with the poset of faces of $\underline{P}$.
The poset of strata of $S$ is therefore identified with the poset of faces of $\underline{P}$ which are not 
contained in the union of facets $\bigcup (x_i=0) \subset \Delta(r,n)$ corresponding to $\bigcup B_{i,T}^- \subset T$.
\end{rem}

\begin{rem} Let $S'$ be a component of $S$ and $B'$ the divisor on $S'$ given by the restriction of 
$B$ and the double locus. Then, by Section~\ref{lcmodel}, $(S',B')$ is the log canonical model of the complement of a 
hyperplane arrangement.
\end{rem}

Let $\omega_p$ denote the relative dualising sheaf of $p : \bS \rightarrow M$.

\begin{thm}\label{dualisingsheaf}
$\omega_p(\bB)$ is the restriction of the Pl\"ucker line bundle on $G_e \times M$.
\end{thm}

\begin{lem}\label{toricdualisingsheaf}
Let $T \curvearrowright X/S$ be a flat family of reduced stable toric varieties of dimension $d$.
Let $B$ be the relative toric boundary of $X/S$ and $M=\Hom(T,\bG_m)$.
There is a canonical isomorphism $\omega_{X/S} \cong  \cO_X(-B) \otimes \wedge^{d}M$.
\end{lem}

\begin{proof}[Proof of Lemma~\ref{toricdualisingsheaf}]
Let $X^0 \subset X$ be the smooth locus of $X/S$.
The torus action induces a map $\Omega_{X^0/S} \rightarrow \cO_{X^0} \otimes_k \Lie(T)^* = \cO_{X^0} \otimes_{\bZ} M$
which extends to an isomorphism $\Omega_{X^0/S}(\log B) \rightarrow \cO_{X^0} \otimes M$ 
(cf. \cite[p.~116, Prop.~3.1]{Oda88}). 
Taking top exterior powers we obtain an isomorphism $\omega_{X^0/S}(B) \rightarrow \cO_{X^0} \otimes \wedge^{d}M$, 
and twisting by $\cO_X(-B)$ an isomorphism $\omega_{X^0/S} \rightarrow \cO_{X^0}(-B) \otimes \wedge^{d}M$.
We claim this extends to an isomorphism $\omega_{X/S} \rightarrow \cO_X(-B) \otimes \wedge^{d}M$.
Since $\omega_{X/S}$ is flat over $S$ and has $S_2$ fibres it satisfies a relative $S_2$ property, namely
$\omega_{X/S}=j_{\star} \omega_{X^1/S}$ for $j : X^1 \subset X$ an open subscheme such that the complement has fibres of 
codimension at least $2$ (see \cite[Lem.~A.3]{Hacking04}). 
Similarly for $\cO_X(-B)$.
So, it is enough to check the claim on the open locus $X^1 \subset X$ given by the complement of the torus orbits 
of codimension at least $2$ in the fibres.
At a point $P \in X^1$, either $X/S$ is smooth, or $P \notin B$ and the fibre is 
\'{e}tale locally isomorphic to $(xy=0) \subset \bA^{d+1}$. 
In the second case, there is a $T$-invariant affine open neighbourhood $U \subset X$ of $P$ 
such that, working \'{e}tale locally on $S$, the family $T \curvearrowright U/S$ is of the form
$$\bG^d_m \curvearrowright ((xy=f) \subset \bA^2_{x,y} \times \bG^{d-1}_m \times S),$$
where $f \in \cO_S$ and the $\bG^d_m$ action on $\bA^2_{x,y} \times \bG^{d-1}_m$ is given by 
$$\bG_m \times \bG^{d-1}_m \ni (t_0,t) : (x,y,t') \mapsto (t_0x,t_0^{-1}y,tt').$$   
We reduce to the case $d=1$, $S=\bA^1_u$, $f=u$, where the result is well known.
\end{proof}

Let $M=\Hom(H,\bG_m)=\sum(x_i=0) \subset \bZ^n$, the characters of $H$, and $N=M^*=\bZ^n/\bZ e$.

\begin{proof}[Proof of Theorem \ref{dualisingsheaf}]
Let $\cU_e$ and $\cQ_e$ denote the universal sub-bundle and quotient bundle on $G_e$, respectively.
We have canonical isomorphisms
$$\omega_p(\bB) \cong \omega_{\bT/M}(\bB) \otimes \wedge^{n-r} \cN_{\bS/\bT} \cong 
\cO_{\bT} \otimes \wedge^{n-1}M \otimes \wedge^{n-r}\cQ_e|_{\bS}$$
by the adjunction formula, Theorem~\ref{transthm}, and Lemma~\ref{toricdualisingsheaf}. The exact sequence
$$0 \rightarrow \cU_e \rightarrow \cO_{G_e} \otimes h \rightarrow \cQ_e \rightarrow 0$$
on $G_e$ yields the isomorphism
$$\cO_{G_e} \otimes \wedge^{n-1}h^* \otimes \wedge^{n-r} \cQ_e  \cong \wedge^{r}\cU_e^* = \cO_{G_e}(1),$$
where $\cO_{G_e}(1)$ is the Pl\"ucker line bundle. Composing with the above isomorphism
using the equality $M \otimes_{\bZ} k =h^*$, we obtain an isomorphism $\omega_p(\bB) \cong \cO_{G_e}(1)|_{\bS}$, 
as required.
\end{proof}

\section{Special sections}

Let $I \subset [n]$ be a subset with $|I|=r-1$ and let $\bB_I$ denote the scheme-theoretic intersection $\bigcap_{i \in I} \bB_i$.

\begin{prop} \label{sections}
$\bB_I \subset \bS$ is a section of $p : \bS \rightarrow M$. 
For each fibre $(S,B)$ of $p$, $S$ is smooth and $B$ has normal crossings at $B_I$. 
\end{prop}

\begin{proof}
Let $(T,B_T)$ be a fibre of $(\bT,\bB_T)/M$ and $(S,B)$ its visible contour.
Write $I=\{i_1,\ldots,i_{r-1}\}$.
The scheme $B_I = \bigcap_{i \in I} B_i \subset S$ is the intersection of the scheme 
$B_{I,T} =\bigcap_{i \in I} B^+_{i,T} \subset T$ with $G_e$. 
The divisor $B^+_{i,T}$ equals the intersection $T \cap G_{e_i}$, where
$G_{e_i} \subset G$ is the locus of subspaces containing $e_i$, by \cite[Prop.~1.6.10]{Kapranov93}.
Thus $$B_{I,T} \subset \bigcap_{i \in I} G_{e_i} = \bP(k^n/ \langle e_i \ | \ i \in I \rangle)= \bP^{\bar{I}}.$$
The subscheme $B_{I,T} \subset T$ corresponds to the face $\Gamma= \bigcap_{i \in I} (x_i=1)$ of $\Delta(r,n)$, 
which equals the $(n-r)$-simplex
$$\conv \{ e_{i_1}+\cdots+e_{i_{r-1}}+e_j \ | \ j \notin I \}.$$
We deduce $B_{I,T}=\bP^{\bar{I}}$ by dimensions.
Hence $B_{I}$ is equal to the point 
$\langle e, e_{i_1},\ldots, e_{i_{r-1}}\rangle \in G(r,n).$
In particular, $\bB_I$ is a section of $p : \bS \rightarrow M$.

Let $\uP$ denote the subdivision of $P=\Delta(r,n)$ associated to $T$.
We show that $T$ is smooth at a general point of $B_{I,T}$ by analysing the subdivision $\underline{P}$ at $\Gamma$. 
The polytope $P$ lies in the affine hyperplane $(\sum x_i = r) \subset \bR^n$, an affine space under 
$M_{\bR}$.
Let $I'=I \cup \{ i_r \}$, some $i_r \notin I$, and fix an embedding $P \subset M_{\bR}$ by identifying the vertex
$e_{i_1}+\cdots+e_{i_r}$ as the origin.
Let $\langle S \rangle$ denote the cone and $\langle S \rangle_{\bR}$ the vector space generated by a set 
$S \subset M_{\bR}$. Consider the quotient cone
$$\sigma := ( \langle P \rangle + \langle \Gamma \rangle_{\bR} ) / 
\langle \Gamma \rangle_{\bR}.$$
We have 
$\langle P \rangle = \langle e_j-e_i \ | \ j \notin I', \ i \in I' \rangle$ and
$\langle \Gamma \rangle = \langle e_j -e_{i_r} \ | \ j \notin I' \rangle.$
So, identifying $M_{\bR}/ \langle \Gamma \rangle_{\bR}$ with
$(x_j=0, j \notin I') \subset M_{\bR}$, we have $$\sigma=\langle e_{i_r} - e_i \ | \ i \in I \rangle.$$ 
In particular, $\sigma$ is simplicial, and the generators of
$\langle P \rangle$ yield a minimal set of generators of $\sigma$.
We claim that there is a unique maximal polytope $P'$ of $\underline{P}$ containing
$\Gamma$. Indeed, the edges of any such $P'$ are also edges of $P$ (since $P'$ is a matroid polytope, see \cite{GS87}),
so the corresponding cone $\sigma' \subset \sigma$ is generated by a collection of edges of $\sigma$.
Hence $\sigma' = \sigma$ because $\sigma$ is simplicial, and $P'$ is unique as claimed.
So $T$ has a unique component $T'$
containing the stratum $B_{I,T}$, and $T'$ is smooth at a general point of $B_{I,T}$
(because $\sigma$  is simplicial and its edges generate the lattice). 
We deduce that $S$ is smooth at $B_I$ by Theorem~\ref{transthm}.
\end{proof}

Recall that $h=k^n/ k \cdot e$, the Lie algebra of $H$.

\begin{thm} \label{residues}
Let $$ \res: p_* \omega_p(\bB) \rightarrow \oplus_{I} p_* \cO_{\bB_I} = \wedge^{r-1} k^n \otimes \cO_M$$
be the canonical map given by taking residues along the special sections.
Let 
$$
c:= \wedge^{r-1} h^* \otimes \cO_M \rightarrow p_* \omega_p(\cB)
$$
be the map defining the embedding $\bS \subset G_e \times M \subset \bP(\wedge^{r-1}h) \times M$.
The composition
$$
\res \circ c: \wedge^{r-1} h^* \otimes \cO_M \rightarrow \wedge^{r-1} k^n \otimes \cO_M
$$
is induced by the inclusion $h^* \subset k^n$, $c$ is
an isomorphism, and $\res$ is an isomorphism onto its image.
\end{thm}

\begin{proof}
Let $I \subset [n]$ be a subset of size $r-1$. Write $I=\{i_1,\ldots,i_{r-1}\}$ where $i_1< \cdots < i_{r-1}$.
The residue map $\omega_p(\bB) \rightarrow \cO_{\bB_I}$ is identified with the restriction of the residue map
$\omega_{\bT/M}(\bB) \otimes \wedge^{n-r} \cQ \rightarrow \cO_{\bB_I}$ on $\bT/M$ via the adjunction
$\omega_p(\bB) \cong \omega_{\bT/M}(\bB) \otimes \wedge^{n-r} \cQ|_{\bS}$.  
We explicitly compute this residue map on $\bT/M$.

Let $\bT^0 \subset \bT$ denote the smooth locus of $\bT/M$ and $\bB_{I,T} = \bigcap_{i \in I} \bB_{i,T}$.
We have $\bB_{I,T}=\bP^{\bar{I}} \times M$ where $\bP^{\bar{I}}=\bP(k^n/ \langle e_i \ | \ i \in I \rangle) \subset G(r,n)$,
see the proof of Proposition~\ref{sections}.
Let $\bB^0_{I,T} \subset \bB_{I,T}$ be the open (relative) toric stratum.
Note $\bB^0_{I,T} \subset \bT^0$ by Proposition~\ref{sections}.
Let $N_I=N / \langle e_{i_1},\ldots,e_{i_{r-1}} \rangle$ and $M_I=N_I^{*} \subset M$.
Thus $N_I \otimes \bG_m$ is the quotient torus acting faithfully on $\bB_{I,T}$. 

The adjunction $\omega_{\bT^0/M}(\bB) \rightarrow \omega_{\bB^0_{I,T}}$ is identified, via the isomorphism of 
Lemma~\ref{toricdualisingsheaf}, with the map
$\cO_{\bT^0} \otimes \wedge^{n-1}M \rightarrow \cO_{\bB^0_{I,T}} \otimes \wedge^{n-r} M_I$ induced by the map
$$\langle e_{i_1} \wedge \cdots \wedge e_{i_{r-1}}, \cdot \rangle : \wedge^{n-1}M \rightarrow \wedge^{n-r}M_I.$$
Indeed, the facet $(x_i=1)$ of $P$ corresponding to $\bB_{i,T}$ has outward normal $e_i \in N$, hence a torus invariant 
differential $d\chi^m/\chi^m$ has residue $\langle e_i, m \rangle$ along $\bB_{i,T}^0 := \bB_{i,T} \cap \bT^0 $.
So, the above map is the Poincar\'{e} residue map for $\bB^0_{I,T} \subset \bT^0$ 
(cf. \cite[p.~120]{Oda88},\cite[p.~87]{Fulton93}).

The section $\bB_I \subset \bB^0_{I,T}$ equals $[e] \times M \subset \bP^{\bar{I}} \times M$,
so $\cQ|_{\bB_I}= N_I \otimes \cO_{\bB_I}$, and $\omega_{\bB^0_{I,T}/M} \cong \cO_{\bB^0_{I,T}} \otimes  \wedge^{n-r}M_I$ 
by Lemma~\ref{toricdualisingsheaf}. 
The residue map $\omega_{\bB^0_{I,T}/M} \otimes \wedge^{n-r}\cQ \rightarrow \cO_{\bB_I}$ is induced by the pairing
$\wedge^{n-r}M_I  \otimes \wedge^{n-r} N_I \rightarrow \bZ$.
We obtain the residue map $\omega_{\bT/M}(\bB) \otimes \wedge^{n-r} \cQ \rightarrow \cO_{\bB_I}$ as the composition
$$\omega_{\bT/M}(\bB) \otimes \wedge^{n-r}\cQ \rightarrow \omega_{\bB^0_T} \otimes \wedge^{n-r}\cQ \rightarrow \cO_{\bB_I}.$$

We deduce that the composition
$$\wedge^{r-1} h^* \otimes \cO_{\bS} \rightarrow \cO_{G_e}(1)|_{\bS} \rightarrow \omega_p(\bB) \rightarrow \cO_{\bB_I}$$
is induced by the map $e_{i_1} \wedge \cdots \wedge e_{i_{r-1}} : \wedge^{r-1}h^* \rightarrow k$.
So, the composition  
$$\wedge^{r-1}h^* \otimes \cO_M \rightarrow p_* \omega_p(\bB) \rightarrow \oplus_{|I|=r-1} p_* \cO_{\bB_I}
= \wedge^{r-1} k^n \otimes \cO_M$$ is induced by the inclusion $\wedge^{r-1}h^* \subset \wedge^{r-1}k^n$ as claimed.
Finally, $p_* \omega_p(\bB)$ is locally free of rank $n-1 \choose r-1$ by Proposition~\ref{vanishing} below, so 
$\wedge^{r-1}h^* \otimes \cO_M \rightarrow p_*\omega_p(\bB)$ is an isomorphism.
\end{proof}

\begin{lem} \label{resolution}
Let $Y$ be a projective stable toric variety.
Let $Y^c$ denote the disjoint union of the strata of $Y$ of codimension $c$ which are not contained in the toric boundary and
$p^c : Y^c \rightarrow Y$ the natural map. There is an exact sequence of $\cO_Y$-modules
\begin{equation}\label{1}
0 \rightarrow \cO_Y \rightarrow p^0_* \cO_{Y^0} \rightarrow p^1_*\cO_{Y^1} \rightarrow \cdots.
\end{equation}
Similarly, let $B^c$  denote the disjoint union of the strata of the toric boundary $B$ of codimension $c$ and
$q^c : B^c \rightarrow B$ the natural map. There is an exact sequence of $\cO_B$-modules
\begin{equation}\label{2}
0 \rightarrow \cO_B \rightarrow q^0_* \cO_{B^0} \rightarrow q^1_*\cO_{B^1} \rightarrow \cdots.
\end{equation}
\end{lem}

\begin{proof}
Let $\uP$ be the subdivision of a lattice polytope $P \subset M_{\bR}$ associated to $Y$, and write $d=\dim Y$.
The sequences are defined as follows.
Fix an orientation of each face $P' \in \uP$. For $P'' \subset P'$ a facet and $Y'' \subset Y'$ the corresponding strata of $Y$,
the map $\cO_{Y'} \rightarrow \cO_{Y''}$ is defined to be the restriction map with sign $+1$ if $P'$ and $P''$ are oriented compatibly 
and $-1$ otherwise.  We assume that each maximal polytope and each boundary facet has the orientation
induced by some fixed orientation of $P$, then the maps $\cO_Y \rightarrow p^0_* \cO_{Y^0}$ and 
$\cO_B \rightarrow q^0_* \cO_{B^0}$ are the restriction maps (no signs).

Let $R$ be the homogeneous coordinate ring of $Y$.
By the definition of stable toric varieties, $R$ is the inverse limit of a system $(R_{\sigma},p_{\sigma\tau})$.
The sequence of homogeneous coordinate rings associated to the sequence (\ref{1}) is the sequence
\begin{equation}\label{3}
0 \rightarrow R \rightarrow R^0 \rightarrow R^1 \rightarrow \cdots
\end{equation}
where $R^c$ is the direct sum of the $R_{\sigma}$ for $\sigma \in \Omega$ an interior cone of codimension $c$, and the maps 
$R_{\sigma} \rightarrow R_{\tau}$ for $\tau \subset \sigma$ a facet are $\pm p_{\sigma \tau}$, with the signs
determined as above. Note that by definition the truncated sequence $0 \rightarrow R \rightarrow R^0 \rightarrow R^1$ is 
exact.

The sequence (\ref{3}) is a direct sum of sequences of $k$-vector spaces
$$0 \rightarrow R_a \rightarrow R^0_a \rightarrow R^1_a \rightarrow \cdots$$
indexed by $a \in \omega \cap A$. Recall that $R_{\sigma,a}= k \cdot \chi^a$ if $a \in \sigma$ and $R_{\sigma,a}=0$ otherwise.
We identify the sequence $R^i_a$ with the complex $C_{d-i}(K,L)$ computing the homology of the pair $(K,L)$ of CW-complexes,
where $K=\underline P$ and $L$ is the subcomplex consisting of polytopes $P' \in \uP$ such that $a \notin \Cone(P')$ 
or $P' \subset \partial P$. Let $\upsilon$ denote the cone of $\Omega$ containing $a$ in its relative interior.
The isomorphism $R^{i}_a \rightarrow C_{d-i}(K,L)$ is given by $$R_{\sigma,a} \ni \chi^a \mapsto a(t_{\sigma\upsilon})[P'],$$
where $\sigma=\Cone(P')$ and $[P']$ denotes the generator of $C_{d-i}(K,L)$ corresponding to $P'$ with its chosen orientation.
(The coefficient $a(t_{\sigma\upsilon}) \in k^{\times}$ ensures that the isomorphism is compatible with the boundary maps). 
For $a \neq 0$, the pair $(K,L)$ is homotopy equivalent to the pair $(B^d, B^d - p)$, where $B^d$ is a ball of dimension $d$ 
and $p \in B^d$ an interior point. So $H_i(K,L)=k$ for $i=d$ and $H_i(K,L)=0$ otherwise. 
Thus the graded piece of the sequence (\ref{3}) of weight $a$ is exact for $a \neq 0$. It follows that 
the sequence (\ref{1}) of sheaves on $Y$ associated to (\ref{3}) is exact.

A similar argument shows that the sequence (\ref{2}) is exact.
Let 
\begin{equation} \label{4}
0 \rightarrow S \rightarrow S^0 \rightarrow S^1 \rightarrow \cdots
\end{equation} 
be the associated sequence of homogeneous coordinate rings.
The sequence $S^i_a$ is identified with $C_{d-1-i}(K,L)$, where $K$ is the subcomplex of $\uP$ with support $\partial P$ and 
$L \subset K$ is the subcomplex of faces $P'$ such that $a \notin \Cone(P')$.
For $a \neq 0$, the pair $(K,L)$ is homotopy equivalent to $(S^{d-1},S^{d-1}-p)$, where $S^{d-1}$ is a sphere of dimension 
$(d-1)$ and $p \in S^{d-1}$ a point.
We deduce that the graded piece of the sequence (\ref{4}) of weight $a$ is exact for $a \neq 0$, and 
the sequence (\ref{2}) of sheaves on $Y$ associated to (\ref{4}) is exact, as required.
\end{proof}

\begin{prop} \label{vanishing}
For each fibre $(S,B)$ of $(\bS,\bB)/M$, $\dim_k H^0(\omega_S(B)) = {n-1 \choose r-1}$ and $H^i(\omega_S(B))=0$ for 
$i > 0$. Thus $p_*\omega_p(\bB)$ is locally free of rank $n-1 \choose r-1$ and commutes with base change.
\end{prop}

\begin{proof}
By Serre duality,
$$H^i(\omega_S(B))=\Ext^{r-1-i}(\omega_S(B),\omega_S)^*=H^{r-1-i}(\cO_S(-B))^*,$$
using $S$ Cohen-Macaulay and $\omega_S(B)$ invertible.
We calculate the cohomology groups $H^i(\cO_S(-B))$ using the exact sequence
$$0 \rightarrow \cO_S(-B) \rightarrow \cO_S \rightarrow \cO_B \rightarrow 0.$$
We compute below that $H^i(\cO_S)=0$ for $i>0$, $H^i(\cO_B)=0$ for $0<i<r-2$ and $\dim_k H^{r-2}(\cO_B) ={n-1 \choose r-1}$,
thus $H^i(\cO_S(-B))=0$ for $i<r-1$ and $\dim_k H^{r-1}(\cO_S(-B))={n-1 \choose r-1}$, as required.

Let $(T,B_T)$ be the fibre of $(\bT,\bB_T)/M$ associated to $(S,B)$.
Let $T^c$ denote the disjoint union of the strata of $T$ of codimension $c$ which are not contained in the boundary $B_T$ and
$p^c : T^c \rightarrow T$ the natural map. By Lemma~\ref{resolution}, there is an exact sequence 
$$0 \rightarrow \cO_T \rightarrow p^0_* \cO_{T^0} \rightarrow p^1_* \cO_{T^1} \rightarrow \cdots.$$
Defining $p^c : S^c \rightarrow S$ analogously, we obtain an exact sequence
$$0 \rightarrow \cO_S \rightarrow p^0_* \cO_{S^0} \rightarrow p^1_* \cO_{S^1} \rightarrow \cdots$$
by restriction, using smoothness of $H \times S \rightarrow T$.
For each stratum $S'$ of $S$ we have $H^i(\cO_{S'})=0$ for $i>0$ by Lemma 4.6. So $H^i(\cO_S)$ is the $i$th cohomology of the complex
$$0 \rightarrow H^0(\cO_{S^0}) \rightarrow H^0(\cO_{S^1}) \rightarrow \cdots.$$
By Theorem~\ref{transthm}, the non-boundary strata of $S$ are in bijection with the non-boundary strata of $T$. 
Let $K=\uP$, the subdivision of $P$ associated to $T$, and let $L \subset K$ be the subcomplex with support 
$\partial P$. Then the complex $H^0(\cO_{S^i})$ is identified with the complex $C_{n-1-i}(K,L)$ computing the homology of the pair 
$(K,L)$ of CW-complexes  (cf. Proof of Lemma~\ref{resolution}). We deduce that $H^i(\cO_S)=0$ for $i > 0$.

Similarly, we obtain an exact sequence
$$0 \rightarrow \cO_B \rightarrow q^0_* \cO_{B^0} \rightarrow q^1_* \cO_{B^1} \rightarrow \cdots$$
where $q^c : B^c \rightarrow B$ are the strata of $B$ of codimension $c$, and
$H^i(\cO_B)$ is the $i$th cohomology of the complex
$$0 \rightarrow H^0(\cO_{B^0}) \rightarrow H^0(\cO_{B^1}) \rightarrow \cdots.$$
The strata of $B$ are in bijection with the strata of $B_T$ which are not contained in $\bigcup B^-_{i,T}$. Here $B^-_{i,T}$
is the component of $B$ corresponding to the facet \mbox{$(x_i=0)$} of $P$.
Let $K$ denote the subcomplex of $\uP$ with support $\partial P$ and let $L \subset K$ be the subcomplex with support $\bigcup (x_i=0)$. 
Then the complex $H^0(\cO_{B^i})$ is identified with the complex $C_{n-2-i}(K,L)$. 
To compute the homology, we may replace $\uP$ by the 
trivial subdivision. There is then an isomorphism of chain complexes 
$C_{\cdot}(K,L) \rightarrow C_{\cdot}(\Delta^{(n-2)}_{[n]}, \Delta^{(n-r)}_{[n]})$, 
where $\Delta_{[n]}$ denotes the simplex with vertices labelled by $[n]$ and $\Delta^{(m)}_{[n]}$ its $m$-skeleton,
which sends the facet $(x_i=1)$ of $P$ to $\Delta_{[n]\backslash \{i\}}$.
We find $\dim_k H^{n-r}(K,L)= {n-1 \choose r-1}$ and $H^i(K,L)=0$ for $i \neq n-2,n-r$.
Explicitly, $H^{n-r}(K,L)$ is the cokernel of the boundary map $C_{n-r+1}(\Delta_{[n]}) \rightarrow C_{n-r}(\Delta_{[n]})$, 
which may be identified with the map 
$$\wedge^{r-2}k^n \rightarrow \wedge^{r-1}k^n, \ v \mapsto e \wedge v.$$
Then $H^{n-r}(K,L)$ is identified with $\wedge^{r-1}h$ where $h=k^n/k \cdot e$.
We deduce that $\dim_k H^{r-2}(\cO_B)={n-1 \choose r-1}$ and $H^i(\cO_B)=0$ for $0<i<r-2$.
\end{proof}

\begin{lem} Let $S'$ be a closed stratum of a fibre $S$ of the visible
contour family $\bS \rightarrow M$. $S'$ is rational with rational singularities. 
\end{lem}

\begin{proof} By Theorem~\ref{transthm}, $S'$ has singularities no worse than those of the corresponding stratum
of $T$ (the corresponding fibre of $\bT \rightarrow M$), which is a normal toric variety
(and in particular has at worst rational singularities) by Corollary~\ref{toricfamily}. 
By \cite[3.1.9]{Kapranov93}, $S'$ is rational --- it compactifies the complement to a
hyperplane arrangement.
\end{proof}

\section{Very stable pairs} \label{verystable}

\begin{defn}
A very stable pair over a $k$-scheme $T$ is a family $q \colon (\cS,\cB) \rightarrow T$ of toroidal pairs, where 
$\cB=\cB_1+\cdots+\cB_n$, satisfying the following conditions:
\begin{enumerate}
\item $\cS,\cB_1,\ldots,\cB_n$ are flat over $T$ and the sheaf $\omega_q(\cB)$ is a line bundle.
\item For each subset $I \subset [n]$ with $|I|=r-1$, $\cB_I:= \bigcap_{i \in I} \cB_i \subset \cS$ is a section of $q$.
For each fibre $(S,B)$ of $q$, $S$ is smooth and $B$ has normal crossings at $B_I$.
\item The residue map $q_*\omega_q(\cB) \rightarrow \oplus_{I} q_*\cO_{\cB_I} = \wedge^{r-1} k^n \otimes \cO_T$
is an isomorphism onto $\wedge^{r-1}h^* \otimes \cO_T \subset \wedge^{r-1} k^n \otimes \cO_T$.
Let $c : \wedge^{r-1} h^* \otimes \cO_T \rightarrow q_*\omega_q(\cB)$ denote its inverse.
\item The line bundle $\omega_q(\cB)$ and the isomorphism $c$ define an embedding $\cS \subset \bP(\wedge^{r-1}h) \times T$ 
which factors through $G(r-1,h) \times T$.
\item Let $\cT$ denote the sweep closure $\overline{H \cS}$ of $\cS$ in $G(r,n) \times T$ and similarly 
let $\cB^+_{i,T}= \overline{H \cB_i}$ for each $i$. Then the affine cone over $\cT/T$ is a $T$-valued point 
of the toric Hilbert scheme $H^h_S$, and $\cB^+_{i,T}$ is the component of the relative toric boundary 
of $\cT/T$ corresponding to the facet $(x_i=1)$ of $\Delta(r,n)$.
\end{enumerate}
\end{defn}

\begin{rem}
For $H \curvearrowright X$ a group acting on a scheme $X$ and $Y \subset X$ a subscheme of $X$,
the sweep closure $\overline{HY}$ is by definition the scheme-theoretic image of the multiplication 
map $H \times Y \rightarrow X$.
For $f : Z \rightarrow X$ a map of schemes, the scheme-theoretic image of $f$ is the closed subscheme of 
$X$ defined by the ideal sheaf $\cI = \ker (\cO_X \rightarrow f_* \cO_Z)$.
\end{rem}

\begin{thm} \label{modularthm} 
$M$ is a fine moduli space of very stable pairs, with universal family the family of visible contours 
$p : (\bS,\bB) \rightarrow M$.
\end{thm}

\begin{proof}
An arbitrary pullback of the visible contour family $(\bS,\bB)/M$ is a family of very stable pairs by
Thm.~\ref{transthm}, Thm.~\ref{dualisingsheaf}, Prop.~\ref{sections}, Thm.~\ref{residues}, Prop.~\ref{vanishing}, and 
Lemma~\ref{sweepclosure} below. 
It remains to check that $(\bS,\bB)/M$ is universal.
Let $(\cS,\cB)/T$ be a family of very stable pairs, and consider the associated visible contour family
$$(\cS',\cB')=(\overline{H\cS},\overline{H\cB}) \cap G_e \times T$$
which is obtained by pullback from $(\bS,\bB)/M$. Consider the closed embedding $\cS \subset \cS'$.
Let $S \subset S'$ be the restriction to a general fibre; we claim $S=S'$.
Since $S$ and $S'$ are reduced and have pure dimension $r-1$, $S$ is a union of irreducible components of $S'$. 
Each component $S'_j$ of $S'$ is of the form $T_j \cap G_e$, where $T_j$ is a component of the stable toric variety 
$T=\overline{HS}$. Let $x_j$ be a point of $S$ in the interior of the toric variety $T_j$. Then $S'_j$ is the only 
irreducible component of $S'$ containing $x_j$, so $S'_j \subset S$. Hence $S=S'$ as claimed. 
We deduce $\cS=\cS'$ by flatness. The same argument shows $\cB_i=\cB_i'$.
\end{proof}

\begin{lem} \label{sweepclosure}
Let $T \rightarrow M$ be a morphism and let $\cS,\cB_i,\cT,\cB^+_{i,T}$ denote the pullbacks of $\bS,\bB_i,\bT,\bB^+_{i,T}$.
The sweep closures $\overline{H \cS}, \overline{H \cB_i}$ are equal to $\cT,\cB^+_{i,T}$.
\end{lem}

\begin{proof}
The map $H \times \cS \rightarrow \cT$ is smooth, with image $\cT^0:= \cT - \bigcup \cB^-_{i,T}$, where 
$\cB^-_{i,T}=\bB^-_{i,T}|_{\cT}$.
Hence $\overline{H \cS} = \overline{\cT^0}$. Since $\cT/T$ is flat with reduced fibres, 
any embedded component of $\cT$ contains a fibre by \cite[23.2]{Matsumura89}. In particular there are no 
embedded components contained in $\cT - \cT^0$, so $\overline{\cT^0}=\cT$.
The same argument proves $\overline{H\cB_i}=\cB^+_{i,T}$.
\end{proof}

\section{Example} \label{excpt}

We show that, for $(r,n)=(3,9)$, $M$ has an irreducible component besides the closure of $M^0$.
Moreover, this component is not contained in the image of the Lafforgue space $\oo$ (see Section~\ref{Lafspace}). 

We describe a stable pair $(S,B)$ which is a limit of generic arrangements of 9 lines in $\bP^2$ such that the 
deformation space $\Def(S,B)$ is reducible. More precisely, $\Def(S,B)$ has two smooth components $D_1$ and $D_2$
such that $D_1$ parametrises locally trivial deformations and $D_2$ contains the smoothings of $(S,B)$.
Let $P=[(S,B)]$ denote the corresponding point of $M$.
We show that the map of germs $(P \in M) \rightarrow \Def(S,B)$ is an isomorphism, and the image of Lafforgue's 
space 
$\oo$ in $M$ maps isomorphically onto the smoothing component $D_2$. 

Let $\bar{S}=\bP^2$ and let $\bar{B}=\bar{B}_1+\cdots+\bar{B}_9$ be an arrangement of $9$ lines in $\bP^2$ as 
follows: for $i=1,2,3$, the lines $\bar{B}_i$, are in general position, $\bar{B}_{i+3}=\bar{B}_i$,
and $\bar{B}_{2i+3}$ is a generic line through $\bar{B}_i \cap \bar{B}_{i+1 \! \mod 3}$.
Let $(\bar{\cS},\bar{\cB})/T$ be a generic one parameter smoothing of the pair $(\bar{S},\bar{B})$.
Let $\cS \rightarrow \bar{\cS}$ be the birational morphism given by first blowing up the points 
$B_i \cap B_{i+1 \! \mod 3}$, $i=1,2,3$, then blowing up the strict transforms of the lines $B_i$, $i=1,2,3$.
Let $\cB$ denote the strict transform of $\bar{\cB}$ and $(S,B)$ the special fibre of $(\cS,\cB)/T$.
Then $\cS$ is smooth and $S+\cB$ is a simple normal crossing divisor. One checks that the line bundle 
$\omega_S(B)=\omega_{\cS/T}(\cB)|_S$ is ample. Thus $(S,B)$ is a stable pair.

The deformation space of the surface $S$ may be computed using the results of \cite{Friedman83}.
We find that $\Def S$ is the union of two smooth curve germs $V_1$ and $V_2$ which intersect transversely. 
Here $V_1$ parametrises locally trivial deformations of $S$, and $V_2$ gives the (essentially unique) 1-parameter 
smoothing. The forgetful map $F: \Def(S,B) \rightarrow \Def(S)$ is smooth since $B_i$ is Cartier and 
$H^1(\cN_{B_i/S})=0$ for each $i$ (here $\cN_{B_i/S}$ denotes the normal bundle of $B_i$ in $S$). 
Thus $\Def(S,B)$ is a union of two smooth components $D_i=F^{-1}(V_i)$, $i = 1,2$, as claimed.

We briefly explain the existence of locally trivial deformations of $S$. If $S$ is a reducible surface with 
simple normal crossing singularities, there is a canonically defined line bundle $\cO_D(-S)$ on the double curve 
$D$ of $S$ 
given by $\cI_{S_1}|_D \otimes \cdots \otimes \cI_{S_l}|_D$, 
where $S_1,\cdots,S_l$ are the irreducible components of $S$,
and $\cI_{S_i}$ denotes the ideal sheaf of $S_i \subset S$. If $S$ admits a 1-parameter smoothing $\cS/T$ such 
that the total space is smooth, then $\cO_D(-S)$ is isomorphic to $\cO_D$ 
(because $\cO_D(-S)=\cO_{\cS}(-S)|_D$ and $\cO_{\cS}(-S) \cong \cO_{\cS}$). If $S'$ is a locally trivial 
deformation of $S$, the line bundle $\cO_{D'}(-S')$ lies in $\Pic^0(D')$ but is nontrivial in general.
In our example, $\Pic^0(D) \cong \bG_m$ 
(because $D$ is a union of rational components and contains a unique cycle),
and there are locally trivial deformations $S'$ of $S$ given by changing 
the glueing of the components of $S$ such that $\cO_{D'}(-S')$ is a nontrivial line bundle on $D'$.

We show that the map $(P \in M) \rightarrow \Def(S,B)$ is an isomorphism. By Theorem~\ref{modularthm}, 
it is a closed embedding, and its image contains the smoothing component $D_2$.
It remains to prove that a general fibre over the component $D_1$ of $\Def(S,B)$ is a fibre of the visible 
contour family $(\bS,\bB)/M$. Let $(S,B)$ be an arbitrary fibre over $D_1$. 
The surface $S$ may be identified with the stable toric variety defined by a subdivision of the 
standard triangle of side length $6$ (see the figure) and some glueing data.
\begin{figure}
\begin{center}
\includegraphics[width=0.5\textwidth]{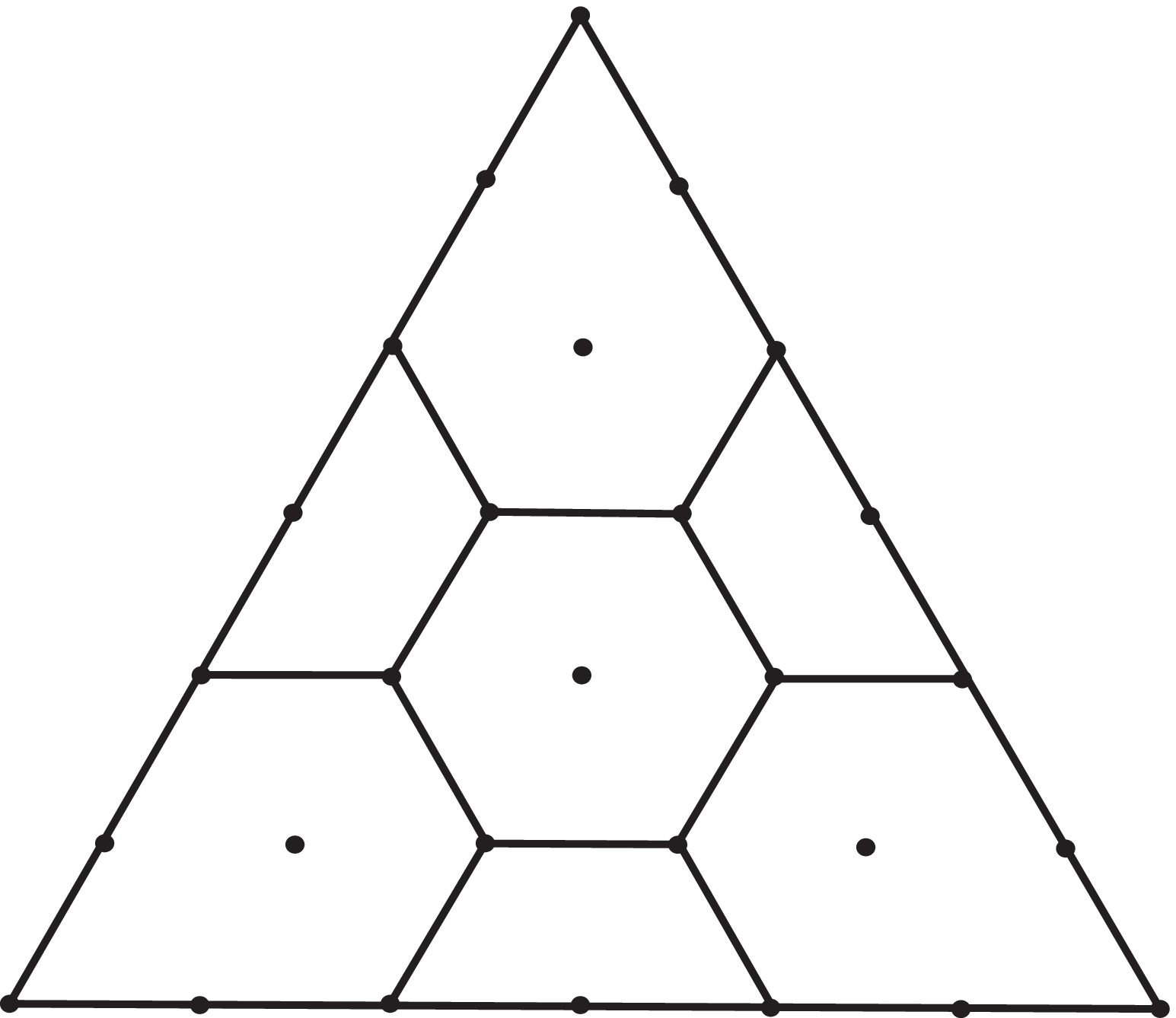}
\end{center}
\end{figure}
The torus action determines a locally free sheaf $\Omega_S(\log)$ on $S$
obtained by glueing the locally free sheaves $\Omega_{S_i}(\log \Delta_i)$ on the components $S_i$ 
at the double locus (here $\Delta_i$ denotes the double locus on $S_i$). There is a natural map
$\Omega_S \rightarrow \Omega_S(\log)$. Let $\Omega_S(\log B)$ be the $\cO_S$-module generated
by $\Omega_S(\log)$ and $\{\frac{df}{f} \ | \ f \in \cO^{\times}_U \}$, where $U=S \setminus B$. Then 
$\Omega_S(\log B)$ is also locally free, and there is an exact sequence
$$0 \rightarrow \Omega_S(\log) \rightarrow \Omega_S(\log B) \rightarrow \oplus \cO_{B_i} 
\rightarrow 0$$
where the last map is given by taking residues along the $B_i$.
The residue map induces an isomorphism 
$H^0(\Omega_S(\log B)) \rightarrow h^* = (\sum x_i =0) \subset k^n$. 
This defines an embedding $(S,B) \subset G(r-1,h)=G(r-1,n-1)_e$. For $S'$ a 
component of $S$, let $B'$ denote the divisor on $S'$ given by the restriction of $B$ and the double locus. 
Then $U'= S' \backslash B'$
is the complement of a hyperplane arrangement, $(S',B')$ is the log canonical model of $U'$, 
and $\Omega_{S'}(\log B') = \Omega_S(\log B)|_{S'}$. 
One checks that the induced map $h^* \rightarrow H^0(\Omega_{S'}(\log B'))$ coincides with the map of 
Theorem~\ref{lcmodelthm}. Thus the locus $\overline{HS'}$ in $G(r,n)$ is the closure of
a single $H$-orbit.
The weight polytopes $P' \subset P=\Delta(r,n)$ of the orbit closures $\overline{HS'}$ define a subdivision of 
$P$ (because this only depends on the combinatorial type of $(S,B)$, and holds for the fibre over $0\ \in D_2$).
Hence $\overline{HS}$ defines a point of the toric Hilbert scheme $H^h_S$, and $(S,B)$ is its visible contour.
Thus $(S,B)$ is a fibre of the visible contour family over $M$, as required.

The Lafforgue space $\oo$ is a moduli space of varieties with log structures.
We refer to \cite{Kato89} for background on log structures.
Given a pair $[(S,B)] \in M$ which lies in the image of $\oo$, a point of $\oo$ over $[(S,B)]$
corresponds to a log structure on $S/k$ which (in particular) determines the divisors $B_i \subset S$. 
In our example, the log structure on $S/k$ is the restriction
of the log structure on the smoothing $\cS/T$ defined by the divisors $S+\cB \subset \cS$ and $0 \in T$.
By \cite{KN94} the log deformations of $S/k$ are parametrised by the component $D_2 \subset \Def(S,B)$, 
thus the germ of $\oo$ at $S/k$ maps isomorphically onto $D_2$.

\medskip
\noindent
Paul Hacking, Department of Mathematics, Yale University, PO Box 208283, New Haven, CT 06520; paul.hacking@yale.edu\\
\\
Sean Keel, Department of Mathematics, University of Texas at Austin, Austin, TX 78712; keel@math.utexas.edu\\
\\
Jenia Tevelev, Department of Mathematics, University of Texas at Austin, Austin, TX 78712; tevelev@math.utexas.edu\\

\end{document}